\documentclass[11pt]{amsart}

\usepackage{amsmath,latexsym,amssymb,amsthm,array,amsfonts}

\usepackage{mathrsfs,dsfont,bbm}

\numberwithin{equation}{section}

\newtheorem{definition}{Definition}[section]
\newtheorem{theorem}{Theorem}[section]
\newtheorem{lemma}{Lemma}[section]

\newtheorem{remark}{Remark}[section]

\newcommand{\E}{{\mathbb E}}

\newcommand\CE{{\mathcal E}}

\def\R{\mathbb{R}}
\newcommand{\W}{\Omega}

\begin{document}
\title{The G-convex Functions Based on the Nonlinear Expectations Defined by G-BSDEs$^*$}

\footnote[0]{${}^{*}$The Project-sponsored by NSFC (11301068),  NSFC
(11171062), NSFC(11371362) and the Fundamental Research Funds for the Central Universities No. 2232014D3-08.}

\author[K. He]{Kun He }


\date{}

\maketitle
\begin{center}
{\footnotesize {\it  hekun\symbol{64}dhu.edu.cn\\
Department of Mathematics\\
 Donghua University\\
2999 North Renmin Rd., Songjiang\\
 Shanghai 201620, P.R. China
}}
\end{center}
\begin{abstract}
In this paper, generalizing the definition of $G$-convex functions defined by Peng \cite{Peng2010} during the construction of $G$-expectations and related properties, we define a group of $G$-convex functions based on the Backward Stochastic Differential Equations driven by $G$-Brownian motions.
\end{abstract}
\section{Introduction}
As we all know, Jensen's inequality is an important result in the theory of linear expectations.
For nonlinear expectations, Jiang \cite{CKJ1,CKJ2,J} talked about a type of $g-$expectation and found out that if this group of
 $g-$expectations satisfy the related Jensen's inequality, then the corresponding generator $g$ should satisfies
the propositions of positive homogenous and subadditivity. After this in 2010, Jia and Peng \cite{JP}
defined a new group of functions as $g-$convexity and gave a necessary and sufficient condition for a $\mathbb{C}^2$ function being
a $g-$convex function.
\par
After the construction of the G-expectations by Peng's work from 2005 to 2010 \cite{Peng2005,Peng2007a,Peng2008a,Peng2008b,Peng2010},
another series of work \cite{STZ,Song1,Song2} aims at solving an opening problem, a G-martingale $M$ which can be decomposed into a sum of a symmetric G-martingale $\bar{M}$ and a decreasing G-martingale $K$, and this problem was solved by \cite{PSZ} in 2012. Then Hu, Ji, Peng and Song defined a new type of
Backward Stochastic Differential Equation driven by G-Brownian motion (G-BSDE) \cite{HJPS1}, proved a related comparison theorem and defined the group of nonlinear expectations by the solutions of G-BSDEs \cite{HJPS2}. Based on their definition of this group of nonlinear expectations and the related comparison theorem of G-BSDEs, He and Hu \cite{HH} proved a representation
theorem for this group of nonlinear expectations and proved some related equivalent conditions between the generator and
related nonlinear expectations. In this paper, we will talk about the $G-$convex function, defined in Peng \cite{Peng2010}, under the framework of the nonlinear expectations defined by the G-BSDEs \cite{HJPS1}. In section \ref{Sec:Preliminary}, we recall some fundamental definitions
and results about G-expectations and G-BSDEs. In section \ref{Sec:Main}, we will prove our main result, giving the equivalent condition of G-convex function under the framework of G-BSDEs.

\section{Preliminary}\label{Sec:Preliminary}

Let us recall some notations for the related spaces of random variables, definitions and results in the construction of G-Brownian motions and G-expectations. The readers may refer to \cite{Peng2007a,Peng2008a,Peng2008b,Peng2010,HJPS1}.
Throughout the paper, for $x\in\R^d$, we denote $|x|=\sqrt{x\cdot x}$ and $\langle x,x\rangle=x\cdot x$.
\begin{definition}\label{def2.1}
Let $\Omega$ be a given set and let $\mathcal{H}$ be a vector
lattice of real valued functions defined on $\Omega$, namely $c\in \mathcal{H}$
for each constant $c$ and $|X|\in \mathcal{H}$ if $X\in \mathcal{H}$.
$\mathcal{H}$ is considered as the space of random variables. A sublinear
expectation $\mathbb{\hat{E}}$ on $\mathcal{H}$ is a functional $\mathbb{\hat
{E}}:\mathcal{H}\rightarrow \mathbb{R}$ satisfying the following properties:
for all $X,Y\in \mathcal{H}$, we have
\item[(a)] Monotonicity: If $X\geq Y$ then $\mathbb{\hat{E}}[X]\geq
\mathbb{\hat{E}}[Y]$;
\item[(b)] Constant preservation: $\mathbb{\hat{E}}[c]=c$;
\item[(c)] Sub-additivity: $\mathbb{\hat{E}}[X+Y]\leq \mathbb{\hat{E}
}[X]+\mathbb{\hat{E}}[Y]$;

\item[(d)] Positive homogeneity: $\mathbb{\hat{E}}[\lambda X]=\lambda
\mathbb{\hat{E}}[X]$ for each $\lambda \geq 0$.
$(\Omega,\mathcal{H},\mathbb{\hat{E}})$ is called a sublinear expectation space.
\end{definition}

\begin{definition}
\label{def2.2} Let $X_{1}$ and $X_{2}$ be two $n$-dimensional random vectors
defined respectively in sublinear expectation spaces $(\Omega_{1}%
,\mathcal{H}_{1},\mathbb{\hat{E}}_{1})$ and $(\Omega_{2},\mathcal{H}%
_{2},\mathbb{\hat{E}}_{2})$. They are identically distributed, denoted
by $X_{1}\overset{d}{=}X_{2}$, if $\mathbb{\hat{E}}_{1}[\varphi(X_{1}%
)]=\mathbb{\hat{E}}_{2}[\varphi(X_{2})]$, for all$\  \varphi \in C_{b.Lip}%
(\mathbb{R}^{n})$, where $C_{b.Lip}(\mathbb{R}^{n})$ denotes the space of
bounded and Lipschitz functions on $\mathbb{R}^{n}$.
\end{definition}

\begin{definition}
\label{def2.3} In a sublinear expectation space $(\Omega,\mathcal{H}%
,\mathbb{\hat{E}})$, a random vector $Y=(Y_{1},\cdot \cdot \cdot,Y_{n})$,
$Y_{i}\in \mathcal{H}$, is said to be independent of another random vector
$X=(X_{1},\cdot \cdot \cdot,X_{m})$, $X_{i}\in \mathcal{H}$ under $\mathbb{\hat
{E}}[\cdot]$, denoted by $Y\bot X$, if for every test function $\varphi \in
C_{b.Lip}(\mathbb{R}^{m}\times \mathbb{R}^{n})$ we have $\mathbb{\hat{E}%
}[\varphi(X,Y)]=\mathbb{\hat{E}}[\mathbb{\hat{E}}[\varphi(x,Y)]_{x=X}]$.
\end{definition}

\begin{definition}
\label{def2.4} ($G$-normal distribution) A $d$-dimensional random vector
$X=(X_{1},\cdot \cdot \cdot,X_{d})$ in a sublinear expectation space
$(\Omega,\mathcal{H},\mathbb{\hat{E}})$ is called $G$-normally distributed if
for each $a,b\geq0$ we have
\[
aX+b\bar{X}\overset{d}{=}\sqrt{a^{2}+b^{2}}X,
\]
where $\bar{X}$ is an independent copy of $X$, i.e., $\bar{X}\overset{d}{=}X$
and $\bar{X}\bot X$. Here the letter $G$ denotes the function
\[
G(A):=\frac{1}{2}\mathbb{\hat{E}}[\langle AX,X\rangle]:\mathbb{S}%
_{d}\rightarrow \mathbb{R},
\]
where $\mathbb{S}_{d}$ denotes the collection of $d\times d$ symmetric matrices.
\end{definition}

Peng \cite{Peng2008b} showed that $X=(X_{1},\cdot \cdot \cdot,X_{d})$ is $G$-normally
distributed if and only if for each $\varphi \in C_{b.Lip}(\mathbb{R}^{d})$,
$u(t,x):=\mathbb{\hat{E}}[\varphi(x+\sqrt{t}X)]$, $(t,x)\in \lbrack
0,\infty)\times \mathbb{R}^{d}$, is the solution of the following $G$-heat
equation:%
\[
\partial_{t}u-G(D_{x}^{2}u)=0,\ u(0,x)=\varphi(x).
\]

The function $G(\cdot):\mathbb{S}_{d}\rightarrow \mathbb{R}$ is a monotonic,
sublinear mapping on $\mathbb{S}_{d}$ and $G(A)=\frac{1}{2}\mathbb{\hat{E}%
}[\langle AX,X\rangle]\leq \frac{1}{2}|A|\mathbb{\hat{E}}[|X|^{2}]$ implies
that there exists a bounded, convex and closed subset $\Gamma \subset
\mathbb{S}_{d}^{+}$ such that
\[
G(A)=\frac{1}{2}\sup_{\gamma \in \Gamma}\mathrm{tr}[\gamma A],
\]
where $\mathbb{S}_{d}^{+}$ denotes the collection of non-negative elements in
$\mathbb{S}_{d}$.

In this paper, we only consider non-degenerate $G$-normal distribution, i.e.,
there exists some $\underline{\sigma}^{2}>0$ such that $G(A)-G(B)\geq
\underline{\sigma}^{2}\mathrm{tr}[A-B]$ for any $A\geq B$.

\begin{definition}
\label{def2.5} i) Let $\Omega=C_{0}^{d}(\mathbb{R}^{+})$ denote the space of
$\mathbb{R}^{d}$-valued continuous functions on $[0,\infty)$ with $\omega
_{0}=0$ and let $B_{t}(\omega)=\omega_{t}$ be the canonical process. Set
\[
L_{ip}(\Omega):=\{ \varphi(B_{t_{1}},...,B_{t_{n}}):n\geq1,t_{1},...,t_{n}%
\in \lbrack0,\infty),\varphi \in C_{b.Lip}(\mathbb{R}^{d\times n})\}.
\]
Let $G:\mathbb{S}_{d}\rightarrow \mathbb{R}$ be a given monotonic and sublinear
function. $G$-expectation is a sublinear expectation defined by
\[
\mathbb{\hat{E}}[X]=\mathbb{\tilde{E}}[\varphi(\sqrt{t_{1}-t_{0}}\xi_{1}%
,\cdot \cdot \cdot,\sqrt{t_{m}-t_{m-1}}\xi_{m})],
\]
for all $X=\varphi(B_{t_{1}}-B_{t_{0}},B_{t_{2}}-B_{t_{1}},\cdot \cdot
\cdot,B_{t_{m}}-B_{t_{m-1}})$, where $\xi_{1},\cdot \cdot \cdot,\xi_{n}$ are
identically distributed $d$-dimensional $G$-normally distributed random
vectors in a sublinear expectation space $(\tilde{\Omega},\tilde{\mathcal{H}%
},\mathbb{\tilde{E}})$ such that $\xi_{i+1}$ is independent of $(\xi_{1}%
,\cdot \cdot \cdot,\xi_{i})$ for every $i=1,\cdot \cdot \cdot,m-1$. The
corresponding canonical process $B_{t}=(B_{t}^{i})_{i=1}^{d}$ is called a
$G$-Brownian motion.

ii) For each fixed $t\in \lbrack0,\infty)$, the conditional $G$-expectation
$\mathbb{\hat{E}}_{t}$ for $\xi=\varphi(B_{t_{1}}-B_{t_{0}},B_{t_{2}}%
-B_{t_{1}},\cdot \cdot \cdot,B_{t_{m}}-B_{t_{m-1}})\in L_{ip}(\Omega)$, without
loss of generality we suppose $t_{i}=t$, is defined by
\[
\mathbb{\hat{E}}_{t}[\varphi(B_{t_{1}}-B_{t_{0}},B_{t_{2}}-B_{t_{1}}%
,\cdot \cdot \cdot,B_{t_{m}}-B_{t_{m-1}})]
\]%
\[
=\psi(B_{t_{1}}-B_{t_{0}},B_{t_{2}}-B_{t_{1}},\cdot \cdot \cdot,B_{t_{i}%
}-B_{t_{i-1}}),
\]
where
\[
\psi(x_{1},\cdot \cdot \cdot,x_{i})=\mathbb{\hat{E}}[\varphi(x_{1},\cdot
\cdot \cdot,x_{i},B_{t_{i+1}}-B_{t_{i}},\cdot \cdot \cdot,B_{t_{m}}-B_{t_{m-1}%
})].
\]

\end{definition}

For each fixed $T>0$, we set%
\[
L_{ip}(\Omega_{T}):=\{ \varphi(B_{t_{1}},...,B_{t_{n}}):n\geq1,t_{1}%
,...,t_{n}\in \lbrack0,T],\varphi \in C_{b.Lip}(\mathbb{R}^{d\times n})\}.
\]
For each $p\geq1$, we denote by $L_{G}^{p}(\Omega)$ (resp. $L_{G}^{p}%
(\Omega_{T})$) the completion of $L_{ip}(\Omega)$ (resp. $L_{ip}(\Omega_{T})$)
under the norm $\Vert \xi \Vert_{p,G}=(\mathbb{\hat{E}}[|\xi|^{p}])^{1/p}$. It
is easy to check that $L_{G}^{q}(\Omega)\subset L_{G}^{p}(\Omega)$ for $1\leq
p\leq q$ and $\mathbb{\hat{E}}_{t}[\cdot]$ can be extended continuously to
$L_{G}^{1}(\Omega)$.

For each fixed $\mathbf{a}\in \mathbb{R}^{d}$, $B_{t}^{\mathbf{a}}%
=\langle \mathbf{a},B_{t}\rangle$ is a $1$-dimensional $G_{\mathbf{a}}%
$-Brownian motion, where $G_{\mathbf{a}}(\alpha)=\frac{1}{2}(\sigma
_{\mathbf{aa}^{T}}^{2}\alpha^{+}-\sigma_{-\mathbf{aa}^{T}}^{2}\alpha^{-})$,
$\sigma_{\mathbf{aa}^{T}}^{2}=2G(\mathbf{aa}^{T})$, $\sigma_{-\mathbf{aa}^{T}%
}^{2}=-2G(-\mathbf{aa}^{T})$. Let $\pi_{t}^{N}=\{t_{0}^{N},\cdots,t_{N}^{N}%
\}$, $N=1,2,\cdots$, be a sequence of partitions of $[0,t]$ such that $\mu
(\pi_{t}^{N})=\max \{|t_{i+1}^{N}-t_{i}^{N}|:i=0,\cdots,N-1\} \rightarrow0$,
the quadratic variation process of $B^{\mathbf{a}}$ is defined by%
\[
\langle B^{\mathbf{a}}\rangle_{t}=\lim_{\mu(\pi_{t}^{N})\rightarrow0}%
\sum_{j=0}^{N-1}(B_{t_{j+1}^{N}}^{\mathbf{a}}-B_{t_{j}^{N}}^{\mathbf{a}}%
)^{2}.
\]
For each fixed $\mathbf{a}$, $\mathbf{\bar{a}}\in \mathbb{R}^{d}$, the mutual
variation process of $B^{\mathbf{a}}$ and $B^{\mathbf{\bar{a}}}$ is defined by%
\[
\langle B^{\mathbf{a}},B^{\mathbf{\bar{a}}}\rangle_{t}=\frac{1}{4}[\langle
B^{\mathbf{a}+\mathbf{\bar{a}}}\rangle_{t}-\langle B^{\mathbf{a}%
-\mathbf{\bar{a}}}\rangle_{t}].
\]

\begin{definition}
\label{def2.6} For fixed $T>0$, let $M_{G}^{0}(0,T)$ be the collection of
processes in the following form: for a given partition $\{t_{0},\cdot
\cdot \cdot,t_{N}\}=\pi_{T}$ of $[0,T]$,
\[
\eta_{t}(\omega)=\sum_{j=0}^{N-1}\xi_{j}I_{[t_{j},t_{j+1})}(t),
\]
where $\xi_{j}\in L_{ip}(\Omega_{t_{j}})$, $j=0,1,2,\cdot \cdot \cdot,N-1$. For
$p\geq1$, we denote by $H_{G}^{p}(0,T)$, $M_{G}^{p}(0,T)$ the completion of
$M_{G}^{0}(0,T)$ under the norms $\Vert \eta \Vert_{H_{G}^{p}}=\{ \mathbb{\hat
{E}}[(\int_{0}^{T}|\eta_{s}|^{2}ds)^{p/2}]\}^{1/p}$, $\Vert \eta \Vert
_{M_{G}^{p}}=\{ \mathbb{\hat{E}}[\int_{0}^{T}|\eta_{s}|^{p}ds]\}^{1/p}$ respectively.
\end{definition}

For each $\eta \in M_{G}^{1}(0,T)$, we can define the integrals $\int_{0}%
^{T}\eta_{t}dt$ and $\int_{0}^{T}\eta_{t}d\langle B^{\mathbf{a}}%
,B^{\mathbf{\bar{a}}}\rangle_{t}$ for each $\mathbf{a}$, $\mathbf{\bar{a}}%
\in \mathbb{R}^{d}$. For each $\eta \in H_{G}^{p}(0,T;\mathbb{R}^{d})$ with
$p\geq1$, we can define It\^{o}'s integral $\int_{0}^{T}\eta_{t}dB_{t}$. In the following  $\langle B\rangle$ denotes the quadratic variation of $B$ (refer to \cite{Peng2010,HJPS1,PSZ}).

Let $S_{G}^{0}(0,T)=\{h(t,B_{t_{1}\wedge t},\cdot \cdot \cdot,B_{t_{n}\wedge
t}):t_{1},\ldots,t_{n}\in \lbrack0,T],h\in C_{b,Lip}(\mathbb{R}^{n+1})\}$. For
$p\geq1$ and $\eta \in S_{G}^{0}(0,T)$, set $\Vert \eta \Vert_{S_{G}^{p}}=\{
\mathbb{\hat{E}}[\sup_{t\in \lbrack0,T]}|\eta_{t}|^{p}]\}^{\frac{1}{p}}$.
Denote by $S_{G}^{p}(0,T)$ the completion of $S_{G}^{0}(0,T)$ under the norm
$\Vert \cdot \Vert_{S_{G}^{p}}$.

We consider the following type of $G$-BSDEs (in this paper we always use
Einstein convention):%
\begin{equation}\label{eq:GBSDE1}
Y_{t}    =\xi+\int_{t}^{T}g(s,Y_{s},Z_{s})ds+\int_{t}^{T}f(s,Y_{s}%
,Z_{s})d\langle B\rangle_{s}   
  -\int_{t}^{T}Z_{s}dB_{s}-(K_{T}-K_{t}),
\end{equation} 
where%

\[
g(t,\omega,y,z),f(t,\omega,y,z):[0,T]\times \Omega_{T}\times
\mathbb{R}\times \mathbb{R}\rightarrow \mathbb{R}%
\]
satisfy the following properties:

\begin{itemize}
\item[(H1)] There exists some $\beta>1$ such that for any $y,z$,
$g(\cdot,\cdot,y,z),f(\cdot,\cdot,y,z)\in M_{G}^{\beta}(0,T)$.

\item[(H2)] There exists some $L>0$ such that
\[
|g(t,y,z)-g(t,y^{\prime},z^{\prime})|+|f(t,y,z)-f(t,y^{\prime},z^{\prime})|\leq L(|y-y^{\prime}|+|z-z^{\prime}|).
\]

\end{itemize}

For simplicity, we denote by $\mathfrak{S}_{G}^{\alpha}(0,T)$ the collection
of processes $(Y,Z,K)$ such that $Y\in S_{G}^{\alpha}(0,T)$, $Z\in
H_{G}^{\alpha}(0,T;\mathbb{R})$, $K$ is a decreasing $G$-martingale with
$K_{0}=0$ and $K_{T}\in L_{G}^{\alpha}(\Omega_{T})$.

\begin{definition}
\label{def3.1} Let $\xi \in L_{G}^{\beta}(\Omega_{T})$, $g$ and $f$ satisfy (H1) and
(H2) for some $\beta>1$. A triplet of processes $(Y,Z,K)$ is called a solution
of equation (\ref{eq:GBSDE1}) if for some $1<\alpha \leq \beta$ the following
properties hold:

\begin{itemize}
\item[(a)] $(Y,Z,K)\in \mathfrak{S}_{G}^{\alpha}(0,T)$;

\item[(b)] $Y_{t}=\xi+\int_{t}^{T}g(s,Y_{s},Z_{s})ds+\int_{t}^{T}%
f(s,Y_{s},Z_{s})d\langle B\rangle_{s}-\int_{t}^{T}Z_{s}%
dB_{s}-(K_{T}-K_{t})$.
\end{itemize}
\end{definition}

\begin{lemma}
\label{the1.1} (\cite{HJPS1}) Assume that $\xi \in L_{G}^{\beta}(\Omega_{T})$
and $g$, $f$ satisfy (H1) and (H2) for some $\beta>1$. Then equation
(\ref{eq:GBSDE1}) has a unique solution $(Y,Z,K)$. Moreover, for any $1<\alpha<\beta$
we have $Y\in S_{G}^{\alpha}(0,T)$, $Z\in H_{G}^{\alpha}(0,T;\mathbb{R})$
and $K_{T}\in L_{G}^{\alpha}(\Omega_{T})$.
\end{lemma}

In this paper, we also need the following assumptions for $G$-BSDE
(\ref{eq:GBSDE1}).

\begin{itemize}
\item[(H3)] For each fixed $(\omega,y,z)\in \Omega_{T}\times \mathbb{R}%
\times \mathbb{R}$, $t\rightarrow g(t,\omega,y,z)$ and $t\rightarrow
f(t,\omega,y,z)$ are continuous.

\item[(H4)] For each fixed $(t,y,z)\in \lbrack0,T)\times \mathbb{R}%
\times \mathbb{R}$, $g(t,y,z)$, $f(t,y,z)\in L_{G}^{\beta}(\Omega
_{t})$ and
\[
\lim_{\varepsilon \rightarrow0+}\frac{1}{\varepsilon}\mathbb{\hat{E}}[\int
_{t}^{t+\varepsilon}(|g(u,y,z)-g(t,y,z)|^{\beta}+
|f(u,y,z)-f(t,y,z)|^{\beta})du]=0.
\]

\item[(H5)] $K_t=\int_0^t\eta_sd\langle B\rangle_s-2\int_0^tG(\eta_s)ds$,
 where $\eta\in M_G^p(0,T)$, $p\geq 1$.

\item[(H6)] For each $(t,\omega,y)\in \lbrack0,T]\times \Omega_{T}%
\times \mathbb{R}$, $g(t,\omega,y,0)=f(t,\omega,y,0)=0$.
\end{itemize}

Assume that $\xi \in L_{G}^{\beta}(\Omega_{T})$, $g$ and $f$ satisfy (H1) and
(H2) for some $\beta>1$. Let $(Y^{T,\xi},Z^{T,\xi},K^{T,\xi})$ be the
solution of $G$-BSDE (\ref{eq:GBSDE1}) corresponding to $\xi$, $g$ and $f$
on $[0,T]$. It is easy to check that $Y^{T,\xi}=Y^{T^{\prime},\xi}$ on $[0,T]$
for $T^{\prime}>T$. Following (\cite{HJPS2}), we  define
a nonlinear expectation as
\[
\mathbb{\mathcal{E}}_{s,t}[\xi]=Y_{s}^{t,\xi}\text{\quad for \quad }0\leq s\leq t\leq T.
\]

\begin{remark}
In \cite{HJPS2,HH} they both define the nonlinear expectation under the assumption(H1), (H2) and (H6), their consistent nonlinear
expectation was defined by $\mathcal{{E}}_{t}[\xi]=Y_{t}^{T,\xi}\text{ for }t\in \lbrack0,T]$. As described by \cite{HJPS2}, under assumption (H6), the nonlinear expectation satisfies, for $T_1<T_2$, $\mathcal{E}_{t,T_1}[\xi]=\mathcal{E}_{t,T_2}[\xi]$. Then the $\mathcal{E}_t[\xi]=\mathcal{E}_{t,T}[\xi]=Y_t^{T,\xi}$ notation is used.
\end{remark}


The classical g-expectations possess many properties that are useful
in finance and economics and became an important risk measure tool in financial mathematics under the complete market case. The nonlinear expectations derived by the G-BSDEs is a useful
generalization of $g$-expectations defined on an incomplete market case.

\section{Main result}\label{Sec:Main}
\begin{definition}\label{def:1}
For $\xi\in \mathbb{L}_G^{\infty}(\Omega_s), ~ 0\leq t\leq s\leq T$, we define the function $h\in C^2(\mathbb{R})$ be $G$-convex, if $\mathcal{E}_{t,s}[h(\xi)]\geq h[\mathcal{E}_{t,s}(\xi)]$.
\end{definition}

\begin{lemma}\label{le:1}(see \cite{HJPS1}).
Let $\xi\in L_T^{\beta}(\W_T)$ and $g, f$ satisfy (H1) and (H2) for some $\beta>1$. Assume that $(Y,Z,K)$ satisfies $(Y,Z)\in S_G^{\alpha}(0,T)\times H_G^{\alpha}(0,T;\R^d)$ and K is a decreasing G-martingale with $K_0=0$ and
$K_T\in L_G^{\alpha}(\W_T)$ for some $1<\alpha<\beta$ is a solution of \eqref{eq:GBSDE1}. Then, there exists a constant $C_{\alpha}>0$ depending on $\alpha, T, G$ and $L$ such that
\begin{equation}\label{eq:Y-esti}
|Y_t|^{\alpha}\leq C_{\alpha}\hat{\E}_t\left[|\xi|^{\alpha}+\left(\int_t^T|h_s^0|ds\right)^{\alpha}\right],
\end{equation}
\begin{equation}\label{eq:Z-esti}
\hat{\E}\left[\left(\int_0^T|Z_s|^2ds\right)^{\alpha/2}\right]\leq C_{\alpha}\left\{\hat{\E}\left[\sup_{t\in[0,T]}|Y_t|^{\alpha}\right]+
\left(\hat{\E}\left[\sup_{t\in[0,T]}|Y_t|^{\alpha}\right]\right)^{1/2}
\left(\hat{\E}\left[\left(\int_0^Th_s^0ds\right)^{\alpha}\right]\right)^{1/2}\right\},
\end{equation}
where $h_s^0=|g(s,0,0)|+|f(s,0,0)|$.
\end{lemma}

\begin{lemma}\label{le:2}(see \cite{HJPS1,Song1})
Let $\alpha\geq 1$ and $\delta>0$ be fixed. Then, there exists a constant C depending on $\alpha$ and $\delta$ such that
\begin{equation}\label{eq:xi-esti}
\hat{\E}\left[\sup_{t\in[0,T]}\hat{E}_t[|\xi|^{\alpha}]\right]\leq C\left\{\left(\hat{\E}\left[|\xi|^{\alpha+\delta}\right]\right)^{\alpha/(\alpha+\delta)}+
\hat{\E}\left[|\xi|^{\alpha+\delta}\right]\right\},
\end{equation}
$\forall \xi\in L_G^{\alpha+\delta}(\W_T)$.
\end{lemma}

Following Theorem 12 in \cite{HH}, we have the following representation.
\begin{lemma}\label{le:3}
Suppose (H1)-(H4) satisfied. Take a polynomial growth function $\Phi\in C^2_b(\mathbb{R})$, $s\in[t,t+\epsilon]$.
$$Y_s=\Phi(B_{t+\epsilon}-B_t)+\int_s^{t+\epsilon}g(r,Y_r,Z_r)dr+\int_s^{t+\epsilon}f(r,Y_r,Z_r)d\langle B\rangle_r-
\int_s^{t+\epsilon}Z_rdB_r-(K_{t+\epsilon}-K_s).$$
Then
\begin{equation}\label{eq:g-represent}
L^2_G-\lim_{\epsilon\rightarrow 0+}\frac{\mathcal{E}_{t,t+\epsilon}[\Phi(B_{t+\epsilon}-B_t)]-\Phi(0)}{\epsilon}
=g(t,\Phi(0),\Phi'(0))+2G(f(t,\Phi(0),\Phi'(0))+\frac12\Phi''(0)).
\end{equation}
\end{lemma}
{\bf Proof: } Let $\tilde{Y}_s=Y_s-\Phi(B_s-B_t)$, then $\tilde{Y}_t=Y_t-\Phi(0)$ and $\tilde{Y}_{t+\epsilon}=
Y_{t+\epsilon}-\Phi(B_{t+\epsilon}-B_t)=0$ hold.

By using Ito's formula,
\begin{equation*}\begin{split}
-d\tilde{Y}_s=&-dY_s+d\Phi(B_s-B_t)\\
=&g(s,Y_s,Z_s)ds+f(s,Y_s,Z_s)d\langle B\rangle_s-Z_sdB_s-dK_s\\
&+\Phi'(B_s-B_t)dB_s+\frac 12\Phi''(B_s-B_t)d\langle B\rangle_s
\end{split}\end{equation*}
\begin{equation*}\begin{split}
\tilde{Y}_s=&0+\int_s^{t+\epsilon}g(r,Y_r,Z_r)dr+\int_s^{t+\epsilon}(f(r,Y_r,Z_r)+\frac12\Phi''(B_r-B_t))d\langle B\rangle_r\\
&-\int_s^{t+\epsilon}(Z_r-\Phi'(B_r-B_t))dB_r-(K_{t+\epsilon}-K_s).
\end{split}\end{equation*}
Let $\tilde{Z}_s=Z_s-\Phi'(B_s-B_t)$ and $\tilde{K}_s=K_s$. Then $(\tilde{Y}_s,\tilde{Z}_s,\tilde{K}_s)$
satisfies the G-BSDE:
\begin{equation*}\begin{split}
\tilde{Y}_s=&0+\int_s^{t+\epsilon}g(r,\tilde{Y}_r+\Phi(B_r-B_t),\tilde{Z}_r+\Phi'(B_r-B_t))dr\\
&+\int_s^{t+\epsilon}(f(r,\tilde{Y}_r+\Phi(B_r-B_t),\tilde{Z}_r+\Phi'(B_r-B_t))+\frac12\Phi''(B_r-B_t))d\langle B\rangle_r\\
&-\int_s^{t+\epsilon}\tilde{Z_r}dB_r-(\tilde{K}_{t+\epsilon}-K_s).
\end{split}\end{equation*}
From Lemma \ref{le:1},
\begin{equation*}\begin{split}
|\tilde{Y}^{\epsilon}_s|^{\alpha}&\leq C_{\alpha}\hat{\E}_s\left[\left(\int_s^{t+\epsilon}(
|g(r,\Phi(B_r-B_t),\Phi'(B_r-B_t))|\right.\right.\\
&\left.\left.+|f(r,\Phi(B_r-B_t),\Phi'(B_r-B_t))|+\frac12|\phi''(B_r-B_t)|
)dr\right)^{\alpha}\right],
\end{split}\end{equation*}
\begin{equation*}\begin{split}
\hat{\E}\left[\left(\int_t^{t+\epsilon}|\tilde{Z}_r^{\epsilon}|^2dr\right)^{\alpha/2}\right]&\leq
C_{\alpha}\left\{\hat{\E}\left[\left(\int_t^{t+\epsilon}(|g(r,\Phi(B_r-B_t),\Phi'(B_r-B_t))|+
\frac12|\Phi''(B_r-B_t)|\right.\right.\right.\\
&\left.\left.\left.+|f(r,\Phi(B_r-B_t),\Phi'(B_r-B_t))|)dr\right)^{\alpha}\right]
+\hat{\E}\left[\sup_{s\in[t,t+\epsilon]}|\tilde{Y}^{\epsilon}_s|^{\alpha}\right]\right\}
\end{split}\end{equation*}
hold for some constant $C_{\alpha}>0$, only depending on $\alpha, T, G$ and $L$.


\begin{multline*}
\int_t^{t+\epsilon}\left(|g(r,0,0)|^{\beta}+|f(r,0,0)|^{\beta}\right)dr\leq  2^{\beta-1}\left\{
\epsilon\left(|g(t,0,0)|^{\beta}+|f(t,0,0)|^{\beta}\right)\right.\\
\left.
+\int_t^{t+\epsilon}\left(
|g(r,0,0)-g(t,0,0)|^{\beta}+|f(r,0,0)-f(t,0,0)|^{\beta}\right)dr\right\}.
\end{multline*}

Together with Lemma \ref{le:2} and assumption (H4), we get
\begin{equation}\label{eq:YZ-esti}
\hat{\E}\left[\sup_{s\in[t,t+\epsilon]}|\tilde{Y}_s^{\epsilon}|^{\alpha}+
\left(\int_t^{t+\epsilon}|\tilde{Z}^{\epsilon}_r|^2\right)^{\alpha/2}\right]\leq C_3\epsilon^{\alpha},
\end{equation}
where $C_3$ depends on $x, y, p,\alpha, \beta, T, G$ and $L$.

Now we prove \eqref{eq:g-represent}. Dividing $\epsilon$, take conditional G-expectations and take limits on both sides of the equation in $L_G^2$ norm, then
$\forall \Phi\in C_b^2(\mathbb{R})$,
\begin{equation*}
\begin{split}
\lim_{\epsilon\rightarrow 0+}\frac{\tilde{Y}_t}{\epsilon}=&\lim_{\epsilon\rightarrow 0+}\frac1{\epsilon}
\hat{E}_t[\tilde{Y}_t+(\tilde{K}_{t+\epsilon}-\tilde{K}_t)]\\
=&\lim_{\epsilon\rightarrow 0+}\frac1{\epsilon}\hat{\E}_t\left[\int_t^{t+\epsilon}g(r,\tilde{Y}_r+\Phi(B_r-B_t),\tilde{Z}_r
+\Phi'(B_r-B_t))dr\right.
\\
&\left.+\int_t^{t+\epsilon}\left(f(r,\tilde{Y}_r+\Phi(B_r-B_t),\tilde{Z}_r+\Phi'(B_r-B_t))
+\frac12\Phi''(B_r-B_t)\right)d\langle B\rangle_r\right]\\
=&\lim_{\epsilon\rightarrow 0+}\frac1{\epsilon}\hat{\E}\left[\int_t^{t+\epsilon}g(r,\Phi(B_r-B_t),\Phi'(B_r-B_t))dr\right.\\
&\left.+\int_t^{t+\epsilon}\left(f(r,\Phi(B_r-B_t),\Phi'(B_r-B_t))+\frac12 \Phi''(B_r-B_t)\right)d\langle B\rangle_r
\right]+L_{\epsilon}
\end{split}\end{equation*}
where
\begin{equation*}\begin{split}
L_{\epsilon}=&\frac1{\epsilon}\left\{\hat{\E}\left[\int_t^{t+\epsilon}
g(r,\tilde{Y}_r+\Phi(B_r-B_t),\tilde{Z}_r+\Phi'(B_r-B_t))dr\right.\right.\\ +&\left.\int_t^{t+\epsilon}\left(f(r,\tilde{Y}_r+\Phi(B_r-B_t),\tilde{Z}_r+\Phi'(B_r-B_t))+
\frac12\Phi''(B_r-B_t)\right)d\langle B\rangle_r\right]\\
-&\hat{\E}\left[\int_t^{t+\epsilon}g(r,\Phi(B_r-B_t),\Phi'(B_r-B_t))dr\right.\\
+&\left.\left.\int_t^{t+\epsilon}\left(f(r,\Phi(B_r-B_t),\Phi'(B_r-B_t))+\frac12\Phi''(B_r-B_t)\right)d\langle B\rangle_r\right]\right\}
\end{split}\end{equation*}

It can be verified that $|L_{\epsilon}|\leq (C_4/\epsilon)\hat{\E}[\int_t^{t+\epsilon}(|\tilde{Y}_r|+|\tilde{Z}_r|)dr]$,
where $C_4$ depends on $G, L$ and $T$. By \eqref{eq:YZ-esti}, we have
\begin{equation*}\begin{split}
\hat{\E}[|L_{\epsilon}|^{\alpha}]\leq&
\frac{C_4^{\alpha}}{\epsilon^{\alpha}}\hat{\E}\left[\left(\int_t^{t+\epsilon}(|\tilde{Y}_r|+
|\tilde{Z}_r|)dr\right)^{\alpha}\right]\\
\leq & \frac{2^{\alpha-1}C_4^{\alpha}}{\epsilon^{\alpha}}\hat{\E}\left[\left(\int_t^{t+\epsilon}|\tilde{Y}_r|dr\right)^{\alpha}+
\left(\int_t^{t+\epsilon}|\tilde{Z}_r|dr\right)^{\alpha}\right]\\
\leq &
2^{\alpha-1}C_4^{\alpha}\left\{\hat{\E}\left[\sup_{s\in[t,t+\epsilon]}|\tilde{Y}_s|^{\alpha}\right]+
\epsilon^{-\alpha/2}\hat{\E}\left[\left(\int_t^{t+\epsilon}|\tilde{Z}_r|^2dr\right)^{\alpha/2}\right]\right\}\\
\leq &
2^{\alpha-1}C_4^{\alpha}C_3(\epsilon^{\alpha}+\epsilon^{\alpha/2}),
\end{split}\end{equation*}
then $L_G^{\alpha}-\lim_{\epsilon\rightarrow 0+}L_{\epsilon}=0$.

We set
\begin{equation*}\begin{split}
M_{\epsilon}=&\frac1{\epsilon}\left\{\hat{\E}_t\left[\int_t^{t+\epsilon}g(r,\Phi(B_r-B_t),\Phi'(B_r-B_t))dr\right.\right.\\
+&\left.\int_t^{t+\epsilon}f(r,\Phi(B_r-B_t),\Phi'(B_r-B_t))+\frac12\Phi''(B_r-B_t)d\langle B\rangle_r\right]\\
-&\left.\hat{\E}_t\left[\int_t^{t+\epsilon}g(r,\Phi(0),\Phi'(0))dr+\int_t^{t+\epsilon}\left(f(r,\Phi(0),\Phi'(0))
+\frac12\Phi''(0)\right)d\langle B\rangle_r
\right]\right\}
\end{split}\end{equation*}

By the Lipschitz condition of function $g$ and $f$, and the polynomial growth of $\Phi\in C_b^2(\R)$,
we have $L_{G}^{\alpha}-\lim_{\epsilon\rightarrow 0+}M_{\epsilon}=0$.

Further we set
\begin{equation*}\begin{split}
N_{\epsilon}=&\frac1{\epsilon}\left\{\hat{\E}_t\left[\int_t^{t+\epsilon}g(r,\Phi(0),\Phi'(0))dr
+\int_t^{t+\epsilon}\left(f(r,\Phi(0),\Phi'(0))+\frac12\Phi''(0)\right)d\langle B\rangle_r\right]\right.\\
&\left.-\hat{\E}_t\left[\int_t^{t+\epsilon}g(r,\Phi(0),\Phi'(0))dr+\int_t^{t+\epsilon}\left(
f(r,\Phi(0),\Phi'(0))+\frac12\Phi''(0)\right)
d\langle B\rangle_r\right]\right\}
\end{split}\end{equation*}
You can check that
\begin{equation*}\begin{split}
|N_{\epsilon}|\leq& (C_7/\epsilon)\hat{\E}_t[\int_t^{t+\epsilon}(|g(r,\Phi(0),\Phi'(0))-
g(t,\Phi(0),\Phi'(0))|+\\
&|f(r,\Phi(0),\Phi'(0))-f(t,\Phi(0),\Phi'(0))|)^{\alpha}dr],
\end{split}\end{equation*}
 where $C_7$ depends on $G$.
Then,
\begin{equation*}\begin{split}
\hat{\E}[|N_{\epsilon}|^{\alpha}]\leq & C_7^{\alpha} \frac1{\epsilon}\hat{\E}\left[\int_t^{t+\epsilon}(
|g(r,\Phi(0),\Phi'(0))-g(t,\Phi(0),\Phi'(0))|\right.\\
&\left.+|f(r,\Phi(0),\Phi'(0))-f(t,\Phi(0),\Phi'(0))|
)^{\alpha}dr\right]\\
\leq & C_7^{\alpha}\left(\frac1{\epsilon}\hat{\E}\left[\int_t^{t+\epsilon}(|g(r,\Phi(0),\Phi'(0))-g(t,\Phi(0),\Phi'(0))|
\right.\right.\\
&\left.\left.+|f(r,\Phi(0),\Phi'(0))-f(t,\Phi(0),\Phi'(0))|)^{\beta}dr\right]\right)^{\alpha/\beta}.
\end{split}\end{equation*}
Take limits from both sides of the above inequality and use assumption (H4), then we have
\begin{equation*}
L_G^{\alpha}-\lim_{\epsilon\rightarrow 0+}N_{\epsilon}=0.
\end{equation*}
At the same time,
\begin{equation*}\begin{split}
\hat{\E}\left[\int_t^{t+\epsilon}g(r,\Phi(0),\Phi'(0))dr+
\int_t^{t+\epsilon}(f(t,\Phi(0),\Phi'(0))+\frac12\Phi''(0))d\langle B\rangle_r\right]\\
=g(t,\Phi(0),\Phi'(0))\epsilon+\hat{\E}_t\left[f(t,\Phi(0),\Phi'(0))\left(\langle B\rangle_{t+\epsilon}-\langle B\rangle_t\right)\right]\\
=\left[g(t,\Phi(0),\Phi'(0))+2G\left((f(t,\Phi(0),\Phi'(0))+\frac12\Phi''(0))\right)\right]\epsilon.
\end{split}\end{equation*}
Then we have
\begin{equation*}\begin{split}
L_G^{\alpha}-\lim_{\epsilon\rightarrow 0+}\frac{\tilde{Y}_t}{\epsilon}=& \lim_{\epsilon\rightarrow 0+}\frac1{\epsilon}\{Y_t-\Phi(0)\}\\
=& g(t,\Phi(0),\Phi'(0))+2G\left(
(f(r,\Phi(0),\Phi'(0)+\frac12\Phi''(0))\right).
\end{split}\end{equation*}
%
The proof is finished.


\begin{theorem}\label{thm:main}
Suppose (H1)-(H4) satisfied. Take a function $h\in C^2$, $\phi\in C_b^2(\R)$ is polynomial growth function and $h(\phi)\in C_b^2(\R)$.
Then $h$ is a G-convex function that is equivalent with 
\begin{multline}\label{eq:G-convex-equivalent}
g(t,h(y),h'(y)z)+2G(f(t,h(y),h'(y)z)+\frac12h''(y)z^2+\frac12h'(y)A)\geq\\
h'(y)g(t,y,z)+2h'(y)G(f(t,y,z)+\frac12A),\quad\textrm{ for all } y,z\in\R \textrm{ and } A\in\R.
\end{multline}
\end{theorem}
{\bf Proof:  Necessary condition: }

Take a function $h\in C^2$ and $\phi\in C^2_b(\mathbb{R})$, with $H(\phi)\in C^2_b(\mathbb{R})$. By Lemma \ref{le:3} we have
\begin{equation}\label{eq:e-h-phi}\begin{split}
L^2_G-& \lim_{\epsilon\rightarrow 0+}\frac{\mathcal{E}_{t,t+\epsilon}[h(\phi(B_{t+\epsilon}-B_t))]-h(\phi(0))}{\epsilon}=
g(t,h(\phi(0)),h'(\phi(0))\phi'(0))\\
&+2G\left(f(t,h(\phi(0)),h'(\phi(0))\phi'(0))+\frac12h''(\phi(0))(\phi'(0))^2+\frac12h'(\phi(0)
)\phi''(0)\right)
\end{split}\end{equation}
and
\begin{equation}\label{eq:phi-lim-rep}\begin{split}
&L_G^2-\lim_{\epsilon\rightarrow 0+}\frac{\CE_{t,t+\epsilon}[\phi(B_{t+\epsilon}-B_t)]-\phi(0)}{\epsilon}\\
&=g(t,\phi(0),\phi'(0))+2G(f(t,\phi(0),\phi'(0))+\frac12\phi''(0)).
\end{split}\end{equation}
Based on \eqref{eq:phi-lim-rep}, we have
\begin{equation}\label{eq:h-e-phi}\begin{split}
&L_G^2-\lim_{\epsilon\rightarrow 0+}\frac{h(\CE_{t,t+\epsilon}[\phi(B_{t+\epsilon}-B_t)])-h(\phi(0))}{\epsilon}\\
&=h'(\phi(0))\left(g(t,\phi(0),\phi'(0))+2G(f(t,\phi(0),\phi'(0))+\frac12\phi''(0))\right).
\end{split}\end{equation}
If $h$ is a G-convex function, from Definition \ref{def:1}, $h$ satisfies $\mathcal{E}_{t,t+\epsilon}[h(\phi(B_{t+\epsilon}-B_t))]\geq h\{\mathcal{E}_{t,t+\epsilon}[\phi(B_{t+\epsilon}-B_t)]\}$.
From \eqref{eq:e-h-phi} and \eqref{eq:h-e-phi} we have
\begin{equation*}\begin{split}
g(t,h(\phi(0)),h'(\phi(0))\phi'(0))&+2G\left(f(t,h(\phi(0)),h'(\phi(0))\phi'(0))+
\frac12h''(\phi(0))(\phi'(0))^2+\frac12h'(\phi(0)
)\phi''(0)\right)\\
&\geq
h'(\phi(0))\left(g(t,\phi(0),\phi'(0))+2G(f(t,\phi(0),\phi'(0))+\frac12\phi''(0))\right)
\end{split}\end{equation*}
Where $(\phi(0),\phi'(0),\phi''(0))$ are arbitrary values in $\mathbb{R}^3$. Then we get \eqref{eq:G-convex-equivalent}.
\par
{\bf Sufficient condition:} 

Following a series of work of Soner, Touzi and Zhang \cite{STZ}, Song \cite{Song1,Song2}, the work of Peng, Song and Zhang \cite{PSZ}
proved a representation theorem of G-martingales in a complete subspace of $L_G^{\alpha}(\W_T)$ $(\alpha\geq 1)$. They proved
the decomposition of G-martingale of $\hat{\E}_t[\xi]$ can be uniquely represented $K_t=\int_0^t\eta_sd\langle B\rangle_s-\int_0^t2G(\eta_s)ds$. And then use the similar Picard approximation approach used in \cite{HJPS1} we can get the corresponding theorem as Theorem \ref{thm:main} for a normal decreasing martingale with $K_0=0$ and $K_T\in \mathbb{L}^{\alpha}(\W_T)$.

Take $\xi\in L^{\infty}_G(\phi(B_t))$,
we need to prove
$$\mathcal{E}_{s,t}[h(\xi)]\geq h[\mathcal{E}_{s,t}(\xi)].$$
Take
$$Y_u=\xi+\int_u^tg(r,Y_r,Z_r)dr+\int_u^tf(r,Y_r,Z_r)d\langle B\rangle_r-\int_u^tZ_rdB_r-(K_t-K_u).$$
Applying Ito's formula,
$$-dh(Y_r)=h'(Y_r)\left[g(r,Y_r,Z_r)dr+f(r,Y_r,Z_r)d\langle B\rangle_r-Z_rdB_r-dK_r\right]-\frac12h''(Y_r)|Z_r|^2d\langle B\rangle_r,$$
then
\begin{equation*}
\begin{split}
h(Y_u)=&h(\xi)+\int_u^th'(Y_r)g(r,Y_r,Z_r)dr
+\int_u^t\left[h'(Y_r)f(r,Y_r,Z_r)-\frac12h''(Y_r)|Z_r|^2\right]d\langle B\rangle_r\\
&-\int_u^th'(Y_r)Z_rdB_r-\int_u^th'(Y_r)dK_r\\
=&h(\xi)+\int_u^tg(r,h(Y_r),h'(Y_r)Z_r)dr+\int_u^tf(r,h(Y_r),h'(Y_r)Z_r)d\langle B\rangle_r-\int_u^th'(Y_r)Z_rdB_r\\
&+\int_u^t\left(h'(Y_r)g(r,Y_r,Z_r)-g(r,h(Y_r),h'(Y_r)Z_r)\right)dr\\
&+\int_u^t\left(h'(Y_r)f(r,Y_r,Z_r)-\frac12h''(Y_r)|Z_r|^2-f(r,h(Y_r),h'(Y_r)Z_r)\right)d\langle B\rangle_r
-\int_u^th'(Y_r)dK_r.
\end{split}
\end{equation*}
Since the decreasing process $K_r$, G-martingale, satisfies (H5), we have
\begin{equation*}\begin{split}
=&h(\xi)+\int_u^tg(r,h(Y_r),h'(Y_r)Z_r)dr+\int_u^tf(r,h(Y_r),h'(Y_r)Z_r)d\langle B\rangle_r-\int_u^th'(Y_r)Z_rdB_r\\
&+\int_u^t\left[h'(Y_r)g(r,Y_r,Z_r)-g(r,h(Y_r),h'(Y_r)Z_r)+2h'(Y_r)G(\eta_r)\right]dr\\
&-\int_u^t\left[-h'(Y_r)f(r,Y_r,Z_r)+\frac12h''(Y_r)|Z_r|^2+f(r,h(Y_r),h'(Y_r)Z_r)+h'(Y_r)\eta_r\right]d\langle B\rangle_r
\end{split}\end{equation*}
\begin{equation*}\begin{split}
=&h(\xi)+\int_u^tg(r,h(Y_r),h'(Y_r)Z_r)dr+\int_u^tf(r,h(Y_r),h'(Y_r)Z_r)d\langle B\rangle_r-\int_u^th'(Y_r)Z_rdB_r\\
&+\int_u^t\left[h'(Y_r)g(r,Y_r,Z_r)-g(r,h(Y_r),h'(Y_r)Z_r)+2h'(Y_r)G(\eta_r)\right.\\
&\left.-2G\left(f(r,h(Y_r),h'(Y_r)Z_r)+
\frac12h''(Y_r)|Z_r|^2+h'(Y_r)(\eta_r-f(r,Y_r,Z_r))\right)\right]dr-(\tilde{K}_t-\tilde{K}_u),
\end{split}\end{equation*}
where
\begin{equation*}\begin{split}
\tilde{K}_t=&-\left\{\int_0^t\left[f(r,h(Y_r),h'(Y_r)Z_r)+\frac12h''(Y_r)|Z_r|^2+h'(Y_r)\left(\eta_r-f(r,Y_r,Z_r)\right)\right]d\langle B\rangle_r\right.\\
&\left.-2\int_0^tG\left[f(r,h(Y_r),h'(Y_r)Z_r)+\frac12h''(Y_r)|Z_r|^2+h'(Y_r)\left(\eta_r-f(r,Y_r,Z_r)\right) \right]  \right\}
\end{split}\end{equation*}
is a decreasing G-martingale.
Denote $\tilde{Y}_u=h(Y_u)$ and $\tilde{Z}_u=h'(Y_u)Z_u$, then
\begin{equation}\label{eq:G-BSDE-G-conv-1}\begin{split}
\tilde{Y}_u=&h(\xi)+\int_u^tg(r,\tilde{Y}_r,\tilde{Z}_r)dr+
\int_u^tf(r,\tilde{Y}_r,\tilde{Z}_r)d\langle B\rangle_r-\int_u^t\tilde{Z}_rdB_r\\
&+\int_u^t\left[h'(Y_r)g(r,Y_r,Z_r)-g(r,\tilde{Y}_r,\tilde{Z}_r)+2h'(Y_r)G(\eta_r)\right.\\
&\left.-2G\left(f(r,\tilde{Y}_r,\tilde{Z}_r)+
\frac12h''(Y_r)|Z_r|^2+h'(Y_r)(\eta_r-f(r,Y_r,Z_r))\right)\right]dr-(\tilde{K}_t-\tilde{K}_u),
\end{split}\end{equation}
we know from the inequality \eqref{eq:G-convex-equivalent} that the fourth integral is less or equal to 0.
\par
On the other hand, $\CE_{s,t}[h(\xi)]$ is the solution of the following G-BSDE,
\begin{equation}\label{eq:G-BSDE-G-conv-2}
\bar{Y}_u=h(\xi)+\int_u^tg(r,\bar{Y}_r,\bar{Z}_r)dr+\int_u^tf(r,\bar{Y}_r,\bar{Z}_r)d\langle B\rangle_r-\int_u^t\bar{Z}_rdB_r
-(\bar{K}_t-\bar{K}_u)
\end{equation}
Applying the comparison theorem of G-BSDEs \cite{HJPS2}, we have
$$\tilde{Y}_u\leq\bar{Y}_u.$$
Since $\tilde{Y}_u=h(Y_u)=h(\CE_{u,t}[\xi])$ and $\bar{Y}_u=\CE_{u,t}[h(\xi)]$, then
$$\CE_{u,t}[h(\xi)]\geq h(\CE_{u,t}[\xi]).$$
\begin{remark}
In this paper, we covered how G-Brownian motion is 1-dimensional case. In fact, the n-dimensional case is also satisfied. The proof does not have any great difference.
\end{remark}

\end{document}